\documentclass[10pt, a4paper,twoside,reqno]{amsart}
\usepackage{a4,amssymb,amsfonts,verbatim,pstricks,pst-plot,xypic,times}
\usepackage{yfonts}

\usepackage{mltex}
\usepackage{amsmath}

\allowdisplaybreaks[2]

\usepackage[english]{babel}
\usepackage{amssymb,latexsym}
\usepackage{amssymb}
\usepackage{latexsym}

\newcommand{\pre}{\Re }

\newcommand{\HH}{{\mathbb H}}
\newcommand{\CC}{{\mathbb C}}

\newcommand{\RR}{{\mathbb R}}
\renewcommand{\Re}{\mathrm{Re}}

\newcommand{\SSS}{{\mathbb S}}

\renewcommand{\phi}{\varphi}

\newcommand{\Spin}{\mathrm{spin}}

\renewcommand{\Re}{\mathrm{Re}}

\newcommand{\tr}{\mathrm{tr}}
\newcommand{\<}{\left\langle}       
\renewcommand{\>}{\right\rangle}

\newcommand{\Spinc}{\mathrm{Spin^c}}

\newcommand{\id}{\mathrm{Id}}

\newcommand{\wM}{\widetilde M}

\newcommand{\lgra}{\longrightarrow}

\newtheorem{example}{Examples}[section]
\newtheorem{thm}{Theorem}[section]
\newtheorem{lemma}[thm]{Lemma}
\newtheorem{prop}[thm]{Proposition}
\newtheorem{cor}[thm]{Corollary}
\newtheorem{remark}[thm]{Remark}
\newtheorem{remarks}[thm]{Remarks}
\newtheorem{definition}[thm]{Definition}
\newtheorem{notation}[thm]{Notation}
\newtheorem{exabout:ample}[thm]{Example}

\def\id{\mathrm{id}}

\def\tr{\mathrm{tr}}
\def\<#1,#2>{\langle\,#1,\,#2\,\rangle}

\def\s#1{\mathrm{\Sigma}_{#1} \mathrm{M}}

\def\epsilon{{\varepsilon}}

\def\phi{{\varphi}}

 \def\tr{\mathrm{tr}}      
\DeclareMathAlphabet{\doba}{U}{msb}{m}{n}

\def\Spin{{\mathop{\rm spin}}}
\def\Spinc{{\mathop{\rm Spin}^c}}

\let\<\langle 
\let\>\rangle

\def\s#1{\mathrm{\Sigma}_{#1} \mathrm{M}}

\newcommand{\definedas}{\mathrel{\raise.095ex\hbox{\rm :}\mkern-5.2mu=}}
\begin{document}

\title[]{Complex and Lagrangian surfaces of the complex projective plane via K\"ahlerian Killing Spin$^c$ spinors}



\author[RN]{Roger Nakad}
\author[JR]{Julien Roth}

\address[R. Nakad]{Notre Dame University-Louaiz\'e, Faculty of Natural and Applied Sciences, Department of Mathematics and Statistics, P.O. Box 72, Zouk Mikael, Lebanon,}
\email{rnakad@ndu.edu.lb}

\address[J. Roth]{LAMA, Universit\'e Paris-Est Marne-la-Vall\'ee, Cit\'e Descartes, Champs sur Marne, 77454 Marne-la-Vall\'ee  cedex 2, France}
\email{julien.roth@u-pem.fr}

%\address[P. Romon]{LAMA, Universit\'e Paris-Est Marne-la-Vall\'ee, Cit\'e Descartes, Champs sur Marne, 77454 Marne-la-Vall\'ee  cedex 2, France}
%\email{pascal.romon@univ-mlv.fr}
\maketitle
\begin{abstract} The complex projective space $\mathbb C P^2$ of complex dimension $2$ has a $\Spinc$ structure  carrying  K\"ahlerian Killing spinors. The restriction of one of these  K\"ahlerian Killing spinors to a surface $M^2$ characterizes  the isometric immersion of $M^2$ into $\mathbb C P^2$ if the immersion is either Lagrangian or complex. 
\end{abstract}

\section{Introduction}

A classical problem in Riemannian geometry is to know when a
Riemannian manifold $(M^n,g)$ can be
isometrically immersed into a fixed Riemannian manifold $(\widetilde {M}^{n+p},\widetilde {g})$. The case of space forms $\RR^{n+1}$, $\SSS^{n+1}$ and
$\HH^{n+1}$ is well-known. In fact, the Gauss, Codazzi and Ricci equations
are necessary and sufficient conditions. In other ambient spaces, the Gauss, Codazzi and Ricci equations
are necessary but not sufficient in general. Some additional conditions may be required like for the case of complex space forms, products, warped products or 3-dimensional homogeneous spaces  (see \cite{Da,Da2,Kow,LO,PT,Roth3}).

\indent In low dimensions, especially for surfaces, another necessary
and sufficient condition is now well-known, namely the existence of
a special spinor field called {\it generalized Killing spinor field}
(\cite{Fr,Mo,La,LR}). These results are the geometrical invariant versions of  previous works on the spinorial Weierstrass representation by  R. Kusner and N. Schmidt, B. Konoplechenko, I. Taimanov and many others (see \cite{KS, Ko, Ta}). This representation was expressed by T.
Friedrich \cite{Fr} for surfaces in $\RR^3$ and then extended to
other $3$-dimensonal (pseudo-)Riemannian manifolds
\cite{Mo,Roth3,Roth4,LR2} as well as for hypersurfaces of 4-dimensional space forms and products \cite{LR} or hypersurfaces of 2-dimensional complex space forms by means of ${\rm Spin}^c$ spinors \cite{NR}.

\indent More precisely, the restriction
$\varphi$ of a parallel spinor field on $\RR^{n+1}$ to an oriented
Riemannian hypersurface $M^{n}$ is a solution of the generalized
Killing equation 
\begin{equation}\label{killing} \nabla_X\varphi=-\frac{1}{2}A(X)\cdot\varphi, 
\end{equation}
where
$``\cdot"$ and $\nabla$ are respectively the Clifford
multiplication and the spin connection on $M^{n}$, the tensor $A$ is the
Weingarten tensor of the immersion and $X$ any vector field on $M$. Conversely, T. Friedrich proved in
\cite{Fr} that, in the two dimensional case, if there exists a
generalized Killing spinor field satisfying Equation
\eqref{killing}, where $A$ is an arbitrary field of symmetric
endomorphisms of $TM$, then $A$ satisfies the  fundamental Codazzi and Gauss equations in the theory of embedded hypersurfaces in a Euclidean space and consequently, $A$ is  the
Weingarten tensor of a local isometric immersion of $M$ into $\RR^3$.
Moreover, in this case, the solution $\varphi$ of the generalized
Killing equation is equivalently a solution of the Dirac equation
\begin{equation}\label{dirac} 
D\varphi=H\varphi, 
\end{equation}
where $D$ denotes the Dirac operator on $M$, $|\varphi|$ is constant and $H$ is a real-valued function (which is the mean curvature of the immersion in $\RR^3$).

More recently, this approach was adapted by the second author, P. Bayard and M.A. Lawn in codimension two, namely, for surfaces in Riemannian 4-dimensional real space forms \cite{BLR}, and then generalized in the pseudo-Riemannian setting \cite{bay,BP} as well as for 4-dimensional products \cite{Roth}. As pointed out in \cite{RR}, this approach coincides with the Weierstrass type representation for surfaces in $\RR^4$ introduced by Konopelchenko and Taimanov \cite{Ko,Ta2}.

The aim of the present article is to provide an analogue for the complex projective space $\CC P^2$. The key point is that, contrary to the case of hypersurfaces of $\CC P^2$ which was considered in \cite{NR}, the use of ${\rm Spin}^c$ parallel spinors is not sufficient. Indeed, for both canonical or anti-canonical ${\rm Spin}^c$ structures, the parallel spinors are always in the positive half-part of the spinor bundle. But, as proved in \cite{BLR} or \cite{Roth}, a spinor with non-vanishing positive and negative parts is required to get the integrability condition of an immersion in the desired target space. For this reason, ${\rm Spin}^c$ parallel spinors are not adapted to our problem. Thus, we make use of real K\"ahlerian Killing spinors. Therefore, our argument holds for the complex projective space and not for the complex hyperbolic space, since $\CC\HH^2$ does not carry a real or  imaginary K\"ahlerian Killing spinor.

We will focus on the case of complex and Lagrangian immersions into $\CC P^2$. These two cases and especially the Lagrangian case are of particular interest in the study of surfaces in $\CC P^2$ (see \cite{HR1,HR2,Ur} and references therein for instance). 

First, consider $(M^2,g)$ an oriented Riemannian $\Spinc$ surface and $E$ an oriented $\Spinc$ vector bundle of rank $2$ over $M$ with scalar product $\langle\cdot,\cdot\rangle_E$ and compatible connection $\nabla^E$. We denote by $F^M$ (resp. $F^E$) the curvature form (an imaginary $2$-form on $M$) of the auxiliary line bundle defining the $\Spinc$ structure on $M$ (resp. on the vector bundle $E$).  For a spinor field $\varphi$, we define $\bar \varphi$  by $\bar \varphi = \varphi_+ - \varphi_-$, where $\varphi^+$ and $\varphi^-$ denote the positive and negative half part of $\varphi$ (see Section $2$). They are the projections of $\varphi$ on the eigensubspaces for the eigenvalues $+1$ and $-1$ of the complex volume form. The aim of the paper is to prove the following two results:

The first theorem gives a spinorial characterization of complex immersions of surfaces in the complex projective space $\CC P^2$.
\begin{thm}\label{thm1}
Let $(M^2,g)$ be an oriented Riemannian surface and $E$ an oriented vector bundle of rank $2$ over $M$ with scalar product $\langle\cdot,\cdot\rangle_E$ and compatible connection $\nabla^E$. We denote by $\Sigma=\Sigma M\otimes\Sigma E$ the twisted spinor bundle. Let $B:TM\times TM\lgra E$ be a bilinear symmetric map, $j:TM\lgra TM$ a complex structure on $M$ and $t:E\lgra E$ a complex structure on $E$. Assume moreover that $t(B(X,Y))=B(X,j(Y))$ for all  $X\in \Gamma(TM)$ and consider $\{e_1, e_2\}$ an orthonormal frame of $TM$. Then, the  following two statements are equivalent
\begin{enumerate}

\item There exists a ${\rm Spin}^c$  structure on $\Sigma M\otimes\Sigma E$ whose auxiliary line bundle's curvature is given by $F^{M +E} (e_1, e_2):=F^M(e_1, e_2) + F^E(e_1, e_2) = 0$ and a spinor field $\varphi\in\Gamma(\Sigma M\otimes\Sigma E)$ satisfying for all $X \in \Gamma(TM)$,
\begin{eqnarray}\label{partspinorcomp}
\nabla_X\varphi&=&-\frac 12\eta(X)\cdot\varphi-\frac{1}{2}X\cdot\varphi+\frac{i}{2}j(X)\cdot\overline{\varphi},
\end{eqnarray}
such that $\varphi^+$ and $\varphi^-$ never vanish and where $\eta$ is given by
$$
\eta(X)=\sum_{j=1}^2e_j\cdot B(e_j,X).
$$
\item There exists a local isometric {\bf complex} immersion of $(M^2,g)$ into $\CC P^2$ with $E$ as normal bundle and second fundamental form $B$ such that the complex structure of $\CC P^2$ over $M$ is given by $j$ and $t$ (in the sense of Proposition \ref{propimmersioncp2}).
\end{enumerate}
\end{thm}

The second theorem is the analogue of Theorem \ref{thm1} for Lagrangian surfaces in $\CC P^2$.
\begin{thm}\label{thm2}
Let $(M^2,g)$ be an oriented Riemannian surface and $E$ an oriented vector bundle of rank $2$ over $M$ with scalar product $\langle\cdot,\cdot\rangle_E$ and compatible connection $\nabla^E$. We denote by $\Sigma=\Sigma M\otimes\Sigma E$ the twisted spinor bundle. Let $B:TM\times TM\lgra E$ be a bilinear symmetric map, $h:TM\lgra E$ a bundle map and $s: E\lgra TM$ the dual map of $h$. Assume moreover that $h$ and $s$ are parallel, $h\circ s=-\id_E$ and $A_{h(Y)}X+s(B(X,Y))=0$, for all $X\in \Gamma(TM)$, where $A_{\nu}:TM\lgra TM$ is defined by $g(A_{\nu}X,Y)=\langle B(X,Y),\nu\rangle_E$ for all $X,Y\in \Gamma(TM)$ and $\nu\in E$. Then, the following two statements are equivalent
\begin{enumerate}

\item There exists a ${\rm Spin}^c$ structure on $\Sigma M\otimes\Sigma E$ whose auxiliary line bundle's curvature is given by $F^{M +E} (e_1, e_2)=-2i$ and a spinor field $\varphi\in \Gamma(\Sigma M\otimes\Sigma E)$  satisfying for all $X\in\Gamma(M)$,
\begin{eqnarray}\label{partspinorlagr}
\nabla_X\varphi&=&-\frac 12\eta(X)\cdot\varphi-\frac{1}{2}X\cdot\varphi+\frac{i}{2}h(X)\cdot\overline{\varphi},
\end{eqnarray}
such that $\varphi^+$ and $\varphi^-$ never vanish and where $\eta$ is given by
$$
\eta(X)=\sum_{j=1}^2e_j\cdot B(e_j,X).
$$
\item There exists a local isometric {\bf Lagrangian} immersion of $(M^2,g)$ into $\CC P^2$ with $E$ as normal bundle and second fundamental form $B$ such that over $M$ the complex structure of $\CC P^2$ is given by $h$ and $s$ (in the sense of Proposition \ref{propimmersioncp2}).
\end{enumerate}
\end{thm}

\section{Preliminaries and Notation}
In this section, we briefly review  some basic facts about K\"ahler geometry and $\Spinc$ structures on manifolds and their submanifolds.
For more details we refer to  \cite{am_lectures, ballman, besse, spin, montiel, HMU, hijazi-lecture, bookspin, Ba, ginoux-these, GM,ginoux-book}.

\subsection{Spin$^c$ structures on K\"ahler-Einstein manifolds}
Let $(M^{n}, g)$ be an $n$-dimensional closed Riemannian $\Spinc$  manifold and denote by $\Sigma M$ its complex spinor bundle, 
which has complex rank equal to $2^{[\frac{n}{2}]}$ . The bundle $\Sigma M$ is endowed with a Clifford multiplication denoted by ``$\cdot$'' and a scalar product denoted by $\<\cdot, \cdot\>$.  Given a $\Spinc$ structure on $(M^{n}, g)$, one can check that the
determinant line bundle $\mathrm{det}(\Sigma M)$ has a root $L$ of index $2^{[\frac
{n}{2}]-1}$. This line bundle $L$ over $M$ is called the auxiliary line bundle associated with the
$\Spinc$ structure.
%From a topological point of view, a Riemannian manifold has a $\Spinc$ structure if and only if there exists a complex line bundle $L$ (which %will be the auxiliary line bundle) on $M$ such that 
%$$\omega_2 (M) = [c_1(L)]_{mod\ \ 2},$$ 
%where $\omega_2(M)$ is the second Steifel-Whitney class of $M$ and $c_1(L)$ is the first Chern class of $L$. 
In the particular case, when the auxiliary line bundle can be chosen to be trivial, the manifold is called a $\Spin$ manifold. 
%In this case, we denote by $\Sigma^{'}M$  the spinor bundle or the $\Spin$ bundle. It can be chosen such that the  auxiliary line bundle is %trivial.
%Locally, a $\Spinc$ structure always exists. In fact, the square root of the auxiliary line bundle
%$L$ and $\Sigma^{'}M$ always exist locally. But, $\Sigma M =\Sigma^{'}M \otimes L^{\frac 12} $ is defined globally. This essentially means %that, while
%the spinor bundle and $L^\frac{1}{2}$ may not exist globally, their tensor product (the $\Spinc$ bundle) is
%defined globally.
The connection $\nabla$ on $\Sigma M$ is the twisted connection of the one on the spin bundle (induced 
by the Levi-Civita connection) and a fixed connection $A$ on $L$. 
%The set of  spinors is denoted by $\Gamma (\Sigma M)$.  
The $\Spinc$ Dirac operator $D$ acting on the space of sections of $\Sigma M$ 
is defined locally by 
%multiplication. In local coordinates:
$D =\sum_{j=1}^{n} e_j \cdot \nabla_{e_j},$
where $\{e_j\}_{j=1,\dots, n}$ is a local orthonormal basis of $TM$. 
%$D$ is a first-order elliptic operator and is formally self-adjoint with respect to the $L^2$-scalar product. 
%satisfying for all  spinors $\varphi,\psi$ on $M$ 
%\begin{align}\label{L2-structure_mod_boundary}
%(D\psi, \varphi)= (\psi, D\varphi),
%\end{align}
%where $(., .)$ is the $L^2$-scalar product given by $(\phi, \psi)=\int_M \langle
%\phi, \psi\rangle dv$,  where $dv$  is
%the Riemannian volume form of $M$.  
%A useful tool when examining the $\Spinc$ Dirac operator is the
%Schr\"{o}dinger-Lichnerowicz formula:
%\begin{eqnarray}
%D^2 = \nabla^*\nabla + \frac 14 S +
%\frac{1}{2}F_A \cdot,
%\label{sl}
%\end{eqnarray}
%where $\nabla^*$ is the adjoint of $\nabla$ with respect to the $L^2$-scalar
%product  and $F_A$ is the curvature (imaginary-valued) $2$-form on $M$ associated to the connection $A$ defined on the auxiliary line %bundle $L$, which acts on spinors by the extension of the Clifford multiplication to differential forms. \\

We recall that the complex volume element $\omega_{\mathbb{C}}=i^{[\frac{n+1}{2}]} e_1\wedge \ldots \wedge e_n$ acts as the identity on the spinor bundle if $n$ is odd. If $n$ is even,  $\omega_{\mathbb C}^2=1$. Thus, under the action of the complex volume element, the spinor bundle  decomposes into the  eigenspaces $\Sigma^{\pm} M$ corresponding to the $\pm 1$ eigenspaces, the {\it positive} and the {\it negative} spinors.

Every K\"ahler manifold $(M^{2m},g,J)$ has a canonical 
$\Spinc$ structure induced by the complex structure $J$. 
 The complexified tangent bundle decomposes into
$T^\mathbb C M = T_{1,0} M\oplus T_{0,1} M,$
 the $i$-eigenbundle (resp. $(-i)$-eigenbundle) of the complex linear extension of $J$. 
For any vector field $X$, we denote by $X^{\pm}:=\frac{1}{2}(X\mp iJ(X))$ its component in $T_{1,0} M$, resp. $T_{0,1} M$.
The spinor bundle of the canonical $\Spinc$ structure is defined by 
$$\Sigma M = \Lambda^{0,*} M =\overset{m}{\underset{r=0}{\oplus}} \Lambda^r (T_{0,1}^* M),$$
and its auxiliary line bundle is  $L = (\mathcal K_M)^{-1}= \Lambda^m (T_{0,1}^* M)$, where $\mathcal K_M=\Lambda^{m,0}M$ is the canonical bundle of $M$. 
The line bundle $L$ has a canonical holomorphic connection, whose curvature form is given by $- i\rho$, 
where $\rho$ is the Ricci form defined, for all vector fields $X$ and $Y$, by $\rho(X, Y) = \mathrm{Ric}(J(X), Y)$ and $\mathrm{Ric}$ denotes the Ricci tensor.  
Similarly, one defines  the so called anti-canonical $\Spinc$ structure, whose spinor bundle is given by 
$\Lambda^{*, 0} M =\oplus_{r=0}^m \Lambda^r (T_{1, 0}^* M)$ and the auxiliary line bundle by $\mathcal K_M$. The spinor bundle of any other $\Spinc$ structure on $M$ can be written as:
$$\Sigma M =  \Lambda^{0, *} M \otimes \mathbb L,$$
where $\mathbb L^2 = \mathcal K_M\otimes L$ and  $L$ is the auxiliary line bundle associated with this $\Spinc$ structure. 
The K\"ahler form $\Omega$, defined as $\Omega(X,Y)=g(J(X),Y)$, acts on $\Sigma M$ via Clifford multiplication.
%and this action is locally  given by:
%$$\Omega\cdot \psi = \frac{1}{2} \sum_{j=1}^{2m} e_j\cdot Je_j\cdot\psi,$$
%for all $\psi\in\Gamma(\Sigma M)$, where $\{e_1, \dots, e_{2m}\}$ is a local orthonormal basis of $\mathrm{TM}$. 
Under this action, the spinor bundle decomposes as follows: 
\begin{equation}\label{decomp}
\Sigma M =\overset{m}{\underset{r=0}{\oplus}} \Sigma_r M,
\end{equation}
where $\Sigma_r M$ denotes the eigenbundle to the eigenvalue $i(2r-m)$ of $\Omega$, of complex rank $\binom{m}{k}$. 
It is easy to see that $\Sigma_r M \subset \Sigma^+ M$ (resp. $\Sigma_r M \subset \Sigma^-M$) if and only if $r$ is even (resp. $r$ is odd).
Moreover, for any $X \in \Gamma(TM)$ and $\varphi \in \Gamma(\Sigma_r M)$, we have $X^+ \cdot\varphi \in \Gamma(\Sigma_{r+1}M)$ 
and $X^-\cdot \varphi \in \Gamma(\Sigma_{r-1} M)$, with the convention $\s{-1}=\s{m+1}=M\times\{0\}$.
Thus, for any $\Spinc$ structure, we have 
$\Sigma_r M = \Lambda^{0, r} M\otimes \Sigma_0 M.$
Hence, $(\Sigma_0M)^2 = \mathcal K_M \otimes L,$ where $L$ is the auxiliary line bundle associated with the $\Spinc$ structure. For example, when the manifold is spin, we have $(\Sigma_0M)^2 = \mathcal K_M$ \cite{hitchin, kirchberg1}.  For the canonical $\Spinc$ structure, since $L = (\mathcal K_M)^{-1}$, it follows that $\Sigma_0 M$ is trivial. This yields the existence of parallel spinors (the constant functions) lying in $\Sigma_0M$, \emph{cf.} \cite{Moro1}.

%On a K\"ahler manifold $(M,g,J)$  endowed with any $\Spinc$ structure,  a spinor of the form $\varphi_r+\varphi_{r+1}\in \Gamma(\Sigma_r M\oplus \Sigma_{r+1} M)$, for some $0\leq r\leq m$, is called a {\it   K\"ahlerian Killing $\Spinc$ spinor} if  there exists a non-zero real constant $\kappa$, such that the following equations are satisfied, for all vector fields $X$,
 %\begin{equation}\label{KKSSdefinition}
%\begin{cases}
  %\begin{split}
   %\nabla_X\varphi_r&=  \kappa \ X^-\cdot\varphi_{r+1},\\
   %\nabla_X\varphi_{r+1}&=  \kappa \  X^+\cdot\varphi_{r}.\\
 % \end{split}
%\end{cases}
 %\end{equation}
%K\"ahlerian Killing spinors lying in $\Gamma(\Sigma_{m} M\oplus \Sigma_{m+1} M) = \Gamma(\Sigma_m M)$  or in  $\Gamma(\Sigma_{-1} M\oplus \Sigma_{0} M) = \Gamma(\Sigma_0 M)$ are just parallel spinors.
In \cite{HMU}, the authors gave examples of $\Spinc$ structures on compact K\"ahler-Einstein  manifolds of positive scalar curvature, which carry  K\"ahlerian Killing $\Spinc$ spinors lying in $\Sigma_r M\oplus \Sigma_{r+1} M$, for $r\neq \frac{m\pm1}{2}$, in contrast to the spin case, where K\"ahlerian Killing spinors may only exist for $m$ odd in the middle of the decomposition \eqref{decomp}. We briefly describe  these $\Spinc$ structures here. If the first Chern class $c_1(\mathcal K_M)$ of the canonical bundle of
the K\"ahler manifold $M$ is a non-zero cohomology class,  the greatest number $p\in \mathbb N^*$ such that
$\frac 1p c_1 (\mathcal K_M) \in H^2 (M, \mathbb Z),$
is called the {\it Maslov index} of the K\"ahler manifold. One can thus consider a $p$-th root of the canonical bundle $\mathcal K_M$, \emph{i.e.}
 a complex line bundle $\mathcal L$, such that $\mathcal L^ p = \mathcal K_M$. In \cite{HMU}, O.~Hijazi, S.~Montiel and F.~Urbano proved the following:
\begin{thm}[Theorem~14, \cite{HMU}]
Let $M$ be a $2m$-dimensional K\"ahler-Einstein compact manifold with scalar
curvature $4m(m+1)$ and Maslov index index $p \in  \mathbb N^*$. For each $0 \leq r \leq m+1$, there exists on $M$ a $\Spinc$ structure with auxiliary line bundle given by $\mathcal L^q$, where  $q = \frac{p}{m+1} (2r-m-1) \in \mathbb Z$, and carrying a K\"ahlerian Killing spinor $\psi_{r-1} + \psi_r \in \Gamma(\Sigma_{r-1} M \oplus \Sigma_{r} M)$, \emph{i.e.}  for all $X \in \Gamma(TM)$, it satisfies the first order system
\begin{equation*}
\begin{cases}
\begin{split}
\nabla_X \psi_{r} = - X^+ \cdot \psi_{r-1},\\
\nabla_X \psi_{r-1}  = - X^-  \cdot \psi_{r}.
\end{split}
\end{cases}
\end{equation*}
\end{thm} 
For example, if $M$ is the complex projective space $\mathbb C P^m$  of complex dimension $m$ and constant holomorphic sectional curvature $4$, then $p = m+1$ and $\mathcal L$ is just the tautological line bundle. We fix $0 \leq r \leq m+1$ and we endow $\mathbb C P^m$ with the $\Spinc$ structure whose auxiliary line bundle  is given by $\mathcal L^q$ where $q = \frac{p}{m+1} (2r-m-1) = 2r-m-1 \in \mathbb Z$. For this $\Spinc$ structure, the space of  K\"ahlerian Killing spinors   in $\Gamma(\Sigma_{r-1}M \oplus \Sigma_rM)$ has dimension $\binom{m+1}{r}$.  In this example, for $r=0$ (resp. $r = m+1$), we get the canonical (resp. anticanonical) $\Spinc$ structure for which K\"ahlerian Killing spinors are just parallel spinors. 

\subsection{Submanifolds of Spin$^c$ manifolds.} 
\label{sec22}
Let $(M^2,g)$ be an oriented Riemannian surface, with a given $\Spinc$ structure, and $E$ an oriented $\Spinc$ vector bundle (see \cite[Chapter II]{spin}) of rank 2 on $M$ with an Hermitian product $\langle\cdot,\cdot\rangle_E$ and a compatible connection $\nabla^E$. We consider the spinor bundle $\Sigma$ over $M$ twisted by $E$ and defined by
$\Sigma=\Sigma M\otimes\Sigma E,$
where $\Sigma M$ and $\Sigma E$ are the spinor bundles of $M$ and $E$ respectively. We endow $\Sigma$ with the spinorial connection $\nabla$ defined by
$$\nabla=\nabla^{\Sigma M}\otimes Id_{\Sigma E}+Id_{\Sigma M}\otimes\nabla^{\Sigma E},$$
where $\nabla^{\Sigma M}$ and $\nabla^{\Sigma E}$ are respectively the $\Spinc$ connections on $\Sigma M$ and $\Sigma E$.
We also define the Clifford product ``$\cdot$" by
$$\left\{\begin{array}{l} X\cdot\varphi=(X\cdot_{_M}\alpha)\otimes\overline\sigma\quad\text{if}\ X\in\Gamma(TM),\\ 
X\cdot\varphi=\alpha\otimes(X\cdot_{_E}\sigma)\quad\text{if}\ X\in\Gamma(E),
\end{array}
\right.$$
for all $\varphi=\alpha\otimes\sigma\in\Gamma(\Sigma M\otimes\Sigma E),$ where ``$\cdot_{_M}$" and ``$\cdot_{_E}$" denote the Clifford multiplications on $\Sigma M$ and on $\Sigma E$ respectively and where $\overline{\sigma}=\sigma^+-\sigma^-$ under the natural decomposition of $\Sigma E=\Sigma^+E\oplus\Sigma^-E$. Here, $\Sigma^+E$ and $\Sigma^-E$ are the eigensubbundles (for the eigenvalue $1$ and $-1$) of $\Sigma E$ for the action of the normal volume element $\omega_{\perp}=i\nu_1\cdot_{E}\nu_2$, where $\{\nu_1,\nu_2\}$ is a local orthonormal frame of $E$. Note that $\Sigma^+M$ and $\Sigma^-M$ are defined similarly by the action of the tangent volume element $\omega=ie_1\cdot_{M} e_2$, where $\{e_1,e_2\}$ is an orthonormal basis of $TM$. The twisted Dirac operator $D$ on $\Gamma(\Sigma)$ is defined by
$$D\varphi=e_1\cdot\nabla_{e_1}\varphi+e_2\cdot\nabla_{e_2}\varphi.$$
We note that $\Sigma$ is also naturally equipped  with a hermitian scalar product $\langle.,.\rangle$ which is compatible with the connection $\nabla$, and thus also with a compatible real scalar product $\pre \langle.,.\rangle.$ We also note that the Clifford product ``$\cdot$" of vectors belonging to $TM\oplus E$ is antihermitian with respect to this hermitian product. Finally, we stress that the four subbundles $\Sigma^{\pm\pm}:=\Sigma^{\pm}M\otimes\Sigma^{\pm}E$ are orthogonal with respect to the hermitian product. We will also consider $\Sigma^+=\Sigma^{++}\oplus\Sigma^{--}$ and  $\Sigma^-=\Sigma^{+-}\oplus\Sigma^{-+}$. Throughout the paper we will assume that the hermitian product is $\CC-$linear w.r.t. the first entry, and $\CC-$antilinear w.r.t. the second entry.

Let $(\wM^4, \widetilde g)$  be a Riemannian $\Spinc$ manifold and $(M^2, g)$ an oriented surface isometrically immersed into $\wM$. We denote by $NM$ the normal bundle of $M$ into $\wM$. As $M$ is an oriented surface, it is also $\Spinc$. We denote by $F^{\widetilde M}$  (resp. $F^M$) the curvature $2$-form of the auxiliary line bundle $L^{\wM}$  (resp. $L$) associated with the $\Spinc$ structure on $\wM$ (resp. $M$). Since the manifolds $M$ and $\wM$ are  $\Spinc$, there exists a $\Spinc$ structure on the bundle $NM$ whose   auxiliary line bundle $L_N$  is given
 by $L_N := ({L})^{-1} \otimes {L^{\wM}}_{|_M}.$ We denote by $\Sigma N$ the $\Spinc$ bundle of $N M$ and let $\Sigma=\Sigma M\otimes \Sigma N$ the spinor bundle over $M$ twisted by $NM$ constructed as above with the associated connection and Clifford multiplication. It is a classical fact that the spinor bundle of $\wM$ over $M$, $\Sigma\wM_{|M}$ identifies with $\Sigma$. Moreover the connections on each bundle are linked by the spinorial Gauss formula: for any $\varphi\in\Gamma(\Sigma)$ and any $X\in \Gamma(TM)$,
\begin{eqnarray}\label{gaussspin}
\widetilde{\nabla}_X\varphi=\nabla_X\varphi+\frac{1}{2}\sum_{j=1,2}e_j\cdot B(X,e_j) \cdot\varphi
\end{eqnarray}
where $B$ is the second fundamental form, $\widetilde{\nabla}$ is the spinorial connection of $\Sigma \wM$ and $\nabla$ is the spinoral connection of $\Sigma$ defined as above and $\{e_1,e_2\}$ is a local orthonormal frame of $TM$. Here ``$\cdot$" is the Clifford product on $\Sigma\wM$ which identifies with the Clifford mulitplication on $\Sigma$.

\section{Immersed surfaces into the complex projective space}\label{sec3}

In this section, we will give the basic facts about immersed surfaces in the complex projective plane and in particular derive a sequence of necessary and sufficient conditions for the existence of such immersions.
\subsection{Compatibility equations}
 Let $(M^2,g)$ be a Riemannian surface isometrically immersed in the $2$-dimensional complex projective space  of constant holomorphic sectional curvature $4c>0$. We denote by $\nabla$ the Levi-Civita connection of $\Big(M^2,g:= \left\langle ., .\right\rangle$\Big), $\widetilde{g}$ \Big(also denoted by $\left\langle ., .\right\rangle$ without ambiguity\Big)  the Fubini-Study metric of $\CC P^2(4c)$ and $\widetilde{\nabla}$ its Levi-Civita connection. Moreover, we denote by $\nabla^{\perp}$ the normal connection and $R^{\perp}$ the normal curvature. First of all, we recall that for any $X, Y, Z \in \Gamma(TM)$ the curvature tensor of $\CC P^2(4c)$ is given by
\begin{eqnarray}\label{curv}
\widetilde{R}(X,Y,Z,W)&=&c\Bigg[\left\langle X,W\right\rangle\left\langle Y,Z\right\rangle -\left\langle X,Z\right\rangle \left\langle Y,W\right\rangle+\left\langle J(X),W\right\rangle\left\langle J(Y),Z \right\rangle\nonumber\\ &&\ \ \ \ \ \ -\left\langle J(X),Z\right\rangle \left\langle J(Y),W\right\rangle
+2\left\langle X,J(Y)\right\rangle \left\langle J(Z),W\right\rangle\Bigg].
\end{eqnarray}
The complex structure $J$ induces the existence of the following four operators $$j:TM\lgra TM,\ h:TM\lgra NM,\ s:NM\lgra TM\ \text{and}\ t:NM\lgra NM$$ defined for any $X\in \Gamma(TM)$ and $\xi\in \Gamma(NM)$ by
\begin{eqnarray}\label{relationfhst}
J(X)=j(X)+h(X)\quad\text{and}\quad
J(\xi)=s(\xi)+t(\xi).
\end{eqnarray} 
From the fact that $J^2=-\mathrm{Id}$ and $J$ is antisymmetric, we get that $j$ and $t$ are antisymmetric and for any $X\in\Gamma(TM)$ and $\xi\in\Gamma(NM)$, we have
\begin{align}
&j^2(X)=-X-s\circ h(X),& \label{relation1.1}\\
&t^2(\xi)=-\xi-h\circ s(\xi),& \label{relation1.2}\\
\label{relation1.3}
&j\circ s(\xi)+s\circ t(\xi)=0,&\\
\label{relation1.4}
&h\circ j(X)+t\circ h(X)=0,&\\
\label{relation1.5}
&\left\langle h(X),\xi \right\rangle=-\left\langle X,s(\xi) \right\rangle.&
\end{align}
We denote by $B:TM\times TM\lgra NM$  the second fundamental form and by $S_{\xi}$  the Weingarten operator associated with $\xi\in \Gamma (NM)$ and defined by $\widetilde{g}(S_{\xi}X,Y)=\widetilde{g}(B(X,Y),\xi)$ for any vectors $X,Y$ tangent to $M$. From the fact that $J$ is parallel, we have
\begin{align}
&(\nabla_Xj)Y=S_{h(Y)}X+s(B(X,Y)),& \label{relation2.1}\\
&\nabla^{\perp}_X(h(Y))-h(\nabla_XY)=t(B(X,Y))-B(X,j(Y)),& \label{relation2.2}\\
&\nabla^{\perp}(t(\xi))-t(\nabla^{\perp}_X\xi)=-B(s(\xi),X)-h(S_{\xi}X),& \label{relation2.3}\\
&\nabla_X(s(\xi))-s(\nabla^{\perp}_X\xi)=-j(S_{\xi}X)+S_{t(\xi)}X,&\label{relation2.4}
\end{align}
where Finally, from \eqref{curv}, we deduce that the Gauss, Codazzi and Ricci equations are respectively given by
\begin{eqnarray}\label{gauss}
R(X,Y)Z&=&c\bigg[\left\langle Y,Z \right\rangle X-\left\langle X,Z\right\rangle Y+\left\langle j(Y),Z\right\rangle j(X)-\left\langle j(X),Z\right\rangle j(Y)\nonumber \\
&& \ \ \ \ \ \ \ \ +2\left\langle X,j(Y)\right\rangle j(Z)\bigg]+S_{B(Y,Z)}X-S_{B(X,Z)}Y,\nonumber
\end{eqnarray}
\begin{eqnarray}\label{codazzi}
(\nabla_XB)(Y,Z)-(\nabla_YB)(X,Z)&=&c\bigg[ \left\langle j(Y),Z\right\rangle h(X)-\left\langle j(X),Z \right\rangle h(Y)+2\left\langle j(X),Y\right\rangle h(Z)\bigg],\nonumber
\end{eqnarray}
\begin{eqnarray}\label{ricci}
R^{\perp}(X,Y)\xi&=&c\bigg[ \left\langle h(Y),\xi\right\rangle h(X)-\left\langle h(X),\xi\right\rangle h(Y)+2\left\langle j(X),Y\right\rangle t(\xi)\bigg] \nonumber \\ &&+B(S_{\xi}Y,X)-B(S_{\xi}X,Y).\nonumber
\end{eqnarray}
In the local orthonormal frames $\{e_1,e_2\}$ and $\{\nu_1,\nu_2 \}$ and for any $k,l\in\{1,2\}$, we set $B_{kl}=B(e_k,e_l)$, $j_{kl}=g(j(e_k),e_l)$,  $t_{kl}=\widetilde g(t(\nu_k),e_l)$ and $h_{kl}=\widetilde{g}(h(e_k),\nu_l)$. Thus, these last equations become:
\begin{equation}\label{Gaussgen}K_M=c + <B_{22},B_{11}>-|B_{12}|^2+3c\  (j_{12})^2,\end{equation}
\begin{equation}\label{Riccigen}K_N=-<[S_{\nu_1},S_{\nu_2}](e_1),e_2>+c\Big (h_{21}h_{12}-h_{11}h_{22}+2j_{12}t_{12}\Big ),\end{equation}
\begin{equation}\label{Codazzigen}(\nabla_{e_1}B)(e_2,e_k)-(\nabla_{e_2}B)(e_1,e_k)=c\Big (j_{2k} h(e_1)-j_{1k}h(e_2)+2j_{12} h(e_k)\Big).\end{equation}
We want to point out that $B_{kl}$ is a normal vector whereas $j_{kl}$, $t_{kl}$ and $h_{kl}$ are real numbers.\\ 
It is clear that Equations \eqref{relation1.1} to \eqref{Codazzigen} are necessary conditions for surfaces in $\CC P^2(4c)$. Conversely, given $(M^2,g)$ a Riemannian surface, $E$ a $2$-dimensional vector bundle over $M$ endowed with a scalar product $\overline{g}$ and a compatible connection $\nabla^{E}$. Let $j:TM\lgra TM,\ h:TM\lgra E,\ s:E\lgra TM\ \text{and}\ t:E\lgra E$ be four tensors. Note that the metric $\widetilde{g}:=\<., .\>$ is defined on $TM\oplus E$ by 
$$\left\{\begin{array}{l}
\widetilde{g}(X,Y)=g(X,Y)\quad \text{for any}\ X,Y\in \Gamma(TM), \\
\widetilde{g}(\nu,\xi)=\overline{g}(\xi,\nu)\quad \text{for any}\ \nu,\xi\in \Gamma(E),\\
\widetilde{g}(X,\nu)=0\quad\text{for any}\ X\in \Gamma(TM)\ \text{and}\ \nu\in \Gamma(E).
\end{array}
\right.$$

\begin{definition}
We say that $(M,g,E,\overline{g},\nabla^E,B,j,h,s,t)$ satisfies the compatibility equations for $\CC P^2(4c)$ if $j$ and $t$ are antisymmetric, the Gauss, the Codazzi and Ricci equations \eqref{Gaussgen} \eqref{Riccigen} \eqref{Codazzigen} and equations \eqref{relation1.1}-\eqref{relation2.4} are fulfilled.
\end{definition}

Now, we can state the following classical {\it Fundamental Theorem} for surfaces of $\CC P^2$, which can be found for instance as a special case of the general result of P. Piccione and D. Tausk \cite[Theorem 8.1 and Example 8.2]{PT}.
\begin{prop}\label{propimmersioncp2}
If $(M,g,E,\overline{g},\nabla^E,B,j,h,s,t)$ satisfies the compatibility equations for $\CC P^2(4c)$ then, there exists an isometric immersion $\Phi:M\lgra\CC P^2(4c)$ such that the normal bundle of $M$ for this immersion is isomorphic to $E$ and such that the second fundamental form $II$ and the normal connection $\nabla^{\perp}$  are given by $B$ and $\nabla^E$. Precisely, there exists a vector bundle isometry $\widetilde{\Phi}:E\lgra NM$ so that
$$II=\widetilde{\Phi}\circ B,$$
$$\nabla^{\perp}\widetilde{\Phi}=\widetilde{\Phi}\nabla^E.$$
Moreover, we have 
$$J(\Phi_*(X))=\Phi_*(j(X))+\widetilde{\Phi}(h(X)),$$
$$J(\widetilde{\Phi}(\xi))=\Phi_*(s(X))+\widetilde{\Phi}(t(\xi)),$$
where $J$ is the canonical complex structure of $\CC P^2(4c)$ and this isometric immersion is unique up to an isometry of $\CC P^2(4c)$.

\end{prop}

\subsection{Special cases}
Two special cases are of particular interest and have been widely studied, the complex and Lagrangian surfaces. 

A surface $(M^2,g)$ of $\CC P^2(4c )$ is said {\bf complex} if the tangent bundle of $M$ is stable by the complex structure of $\CC P^2(4c)$, that is, $J(TM)=TM$. Note that we have necessarily $J(NM)=NM$. Hence, in that case, with the above notations, we have $h=0$, $s=0$ and so $j$ and $t$ are respectively almost complex structures on $TM$ and $NM$. The compatibility equations for complex surfaces of $\CC P^2(4c)$ become
\begin{equation}\label{compatibilitycomplex}
\left\{\begin{array}{l}
h=0,\ s=0,\ j^2=-\id_{TM},\ t^2=-\id_{E}\\
\nabla j=0,\ \nabla^{\perp}t=0\\
t(B(X,Y))-B(X,j(Y))=0,\ \forall X, Y\in \Gamma(TM)\\
K_M=4c+<B_{22},B_{11}>-|B_{12}|^2\\
K_N=-<[S_{\nu_1},S_{\nu_2}](e_1),e_2>+2c\\
(\nabla_{e_1}B)(e_2,e_k)-(\nabla_{e_2}B)(e_1,e_k)=0
\end{array}\right.\end{equation}

A surface $(M^2,g)$ of $\CC P^2(4c)$ is said totally real if $J(TM)$ is transversal to $TM$. In the particular case when $J(TM)=NM$, we say that $(M^2,g)$ is {\bf Lagrangian}. In that case, we have $j=0$ and $t=0$. Hence, the compatibility equations for Lagrangian surfaces of $\CC P^2(4c)$ are
\begin{equation}\label{compatibilitylagrangian}
\left\{\begin{array}{l}
j=0,\ t=0,\ s\circ h=-\id_{TM},\ h\circ s=-\id_{E}\\
\nabla s=0,\ \nabla^{\perp}h=0\\
A_{h(Y)}X+s(B(X,Y))=0,\ \forall X, Y\in TM\\
K_M=c+<B_{22},B_{11}>-|B_{12}|^2\\
K_N=-<[S_{\nu_1},S_{\nu_2}](e_1),e_2>+c(h_{21}h_{12}-h_{11}h_{22})\\
(\nabla_{e_1}B)(e_2,e_k)-(\nabla_{e_2}B)(e_1,e_k)=0
\end{array}\right.\end{equation}

%%%%%%%%%%%%%%%%%%%
%%%%%%%%%%%%%%%%%%%%%
%%%%%%%%%%%%%%%%%%%%%%%
%%%%%%%%%%%%%%%%%%%%%ù
\section{Restriction of a K\"ahlerian Killing Spin$^c$ spinor and the curvature computation}\label{sec4}
%\textcolor{red}{Cette section pourrait peut-\^etre \^etre l\'g\`erement plus d\'etaill\'ee}
We consider a special ${\rm Spin}^c$ structure on $\CC P^2$ of constant holomorphic sectional curvature $4c = 4$ carrying a (real) K\"ahlerian Killing spinor $\varphi$.  For example, on can take $q= -1$ and hence $r=1$. For this structure, the curvature of the auxiliary line bundle is given by $F^{\mathbb CP^2}(X,Y)=-2ig(J(X),Y)$.  There exists a spinor $\varphi = \varphi_0 + \varphi_1 \in \Gamma (\Sigma_0 \mathbb C P^2 \oplus \Sigma_1 \mathbb C P^2)$  satisfies the following:
$$\left\{\begin{array}{l}
\widetilde{\nabla}_X\varphi_0=-X^-\cdot\varphi_1,\\
 \widetilde{\nabla}_X\varphi_1=-X^+\cdot\varphi_0,
\end{array}
\right.$$
Thus, we have
$$\widetilde\nabla_X\varphi=-\frac{1}{2}X\cdot\varphi+\frac{i}{2}J(X)\cdot\overline{\varphi},$$
where $\overline\varphi=\varphi_0-\varphi_1$ is the conjugate of $\varphi$ for the action of the complex volume element $\omega_4^{\CC}=-e_1\wedge e_2\wedge e_3\wedge e_4$. Indeed, $\Sigma_0\mathbb C P^2\subset\Sigma^+\mathbb C P^2$ and $\Sigma_1\mathbb C P^2=\Sigma^-\mathbb C P^2$. Note also that such as spinor is of constant norm and each part $\varphi_0$ and $\varphi_1$ does not have any zeros. Indeed, for instance, if $\varphi_0$ vanishes at one point, then it must vanish everywhere (as it is obtained by parallel transport) and $\varphi_1$ is as a parallel spinor which is not the case for this ${\rm Spin}^c$ structure.

Now, let $M$ be a surface of $\CC P^2$ with normal bundle denoted by $NM$. By the identification of the Clifford multiplications and the $\Spinc$ Gauss formula, we have
$$\nabla_X\varphi=-\frac 12\eta(X)\cdot\varphi-\frac{1}{2}X\cdot\varphi+\frac{i}{2}J(X)\cdot\overline{\varphi}.$$
In intrinsic terms, it can be written as 
\begin{eqnarray}\label{partspinor}
\nabla_X\varphi&=&-\frac 12\eta(X)\cdot\varphi-\frac{1}{2}X\cdot\varphi+\frac{i}{2}j(X)\cdot\overline{\varphi}+\frac{i}{2}h(X)\cdot\overline{\varphi},
\end{eqnarray}
where $\eta$ is given by
\begin{equation}\label{relation eta B}
\eta(X)=\sum_{j=1}^2e_j\cdot B(e_j,X).
\end{equation}
Here $B$ is the second fundamental form of the immersion, and the operators $j$ and $h$ are those introduced in Section \ref{sec3}. After projection on $\Sigma^+$ and $\Sigma^-$ respectively, we deduce immediately that
$$\left\{
\begin{array}{l}
\nabla_X\varphi^+=-\dfrac 12\eta(X)\cdot\varphi^+-\dfrac{1}{2}X\cdot\varphi^--\dfrac{i}{2}j(X)\cdot\varphi^--\dfrac{i}{2}h(X)\cdot\varphi^-,\\ \\
\nabla_X\varphi^-=-\dfrac 12\eta(X)\cdot\varphi^--\dfrac{1}{2}X\cdot\varphi^++\dfrac{i}{2}j(X)\cdot\varphi++-\dfrac{i}{2}h(X)\cdot\varphi^+.
\end{array}
\right.$$
From this, since $\overline{\varphi}=\varphi^+-\varphi^-$, we get
\begin{equation}
\label{nablaphibar}
\nabla_X\overline{\varphi}
=-\frac 12 \eta(X)\cdot\overline{\varphi}+\frac{1}{2}X\cdot\overline{\varphi}-\frac{i}{2}j(X)\cdot\varphi-\frac{i}{2}h(X)\cdot\varphi.
\end{equation}

Now, let us go back to an instrinsic setting by considering $(M^2,g)$ an oriented Riemannian surface and $E$ an oriented vector bundle of rank $2$ over $M$ with scalar product $\langle\cdot,\cdot\rangle_E$ and compatible connection $\nabla^E$. We denote by $\Sigma=\Sigma M\otimes\Sigma E$ the twisted spinor bundle. Let $B:TM\times TM\lgra E$ a bilinear symmetric map and $j:TM\lgra TM,\ h:TM\lgra E$ two tensors. We assume that the spinor field $\varphi\in\Gamma(\Sigma)$ satisfies 
%\textcolor{red}{Peut-\^etre que "Curvature computation" serait mieux comme titre car on n'obtient pas encore les condition d'int\'egrabilit\'e dans cette section, on calcule uniquement la courbure}
 Equation \eqref{partspinor}. We will  compute the spinorial curvature for this spinor field $\varphi$ in the local orthonormal frames $\{e_1,e_2\}$ of $TM$ and $\{\nu_1,\nu_2\}$ of $E$. For a sake of simplicity, we can assume that $\{e_1,e_2\}$ is normal at the point $p\in M$ so that at $p$, $ \nabla e_1=0$, $\nabla e_2=0$ and so $[e_1,e_2]=0$. Hence, at the point $p$, we have
\vspace{-0.2cm}
\begin{eqnarray*}
\nabla_{e_1}\nabla_{e_2}\varphi &=&-\frac{1}{2}\nabla_{e_1}(\eta(e_2))\cdot\varphi+\frac{1}{4}\eta(e_2)\cdot\eta(e_1)\cdot\varphi+\frac{1}{4}\eta(e_2)\cdot e_1\cdot\varphi\\
&&-\frac{i}{4}\eta(e_2)\cdot j(e_1)\cdot\overline{\varphi}-\frac{i}{4}\eta(e_2)\cdot h(e_1)\cdot\overline{\varphi}+\frac{1}{4}e_2\cdot\eta(e_1)\cdot\varphi \\
&&+\frac{1}{4}e_2\cdot e_1\cdot\varphi-\frac{i}{4}e_2\cdot j(e_1)\cdot\overline{\varphi}-\frac{i}{4}e_2\cdot h(e_1)\cdot\overline{\varphi}\\
&&+\frac{i}{2}\nabla_{e_1}(j(e_2))\cdot\overline{\varphi}-\frac{i}{4}j(e_2)\cdot\eta(e_1)\cdot\overline{\varphi}+\frac{i}{4}j(e_2)\cdot e_1\cdot\overline{\varphi}\\
&&
+\frac{1}{4}j(e_2)\cdot j(e_1)\cdot\varphi+\frac{1}{4}j(e_2)\cdot h(e_1)\cdot\varphi+\frac{i}{2}\nabla^{\perp}_{e_1}(h(e_2))\cdot\overline{\varphi}\\
&&-\frac{i}{4}h(e_2)\cdot\eta(e_1)\cdot\overline{\varphi}+\frac{i}{4}h(e_2)\cdot e_1\cdot\overline{\varphi}+\frac{1}{4}h(e_2)\cdot j(e_1)\cdot\varphi\\
&&+\frac{1}{4} h(e_2)\cdot h(e_1)\cdot\varphi\\
&=&
\underbrace{- \frac 12 \nabla_{e_1} (\eta(e_2))\cdot \varphi + \frac 14 \eta(e_2) \cdot\eta(e_1)\cdot \varphi + \frac 14 e_2 \cdot e_1 \cdot\varphi   }_{\mathcal A^{12}_1}\\ &&\underbrace{+ \frac 14  \eta(e_2) \cdot e_1 \cdot\varphi + \frac 14 e_2 \cdot\eta(e_1)\cdot \varphi }_{\mathcal A^{12}_2}\\
&& \underbrace{+ \frac i2 \nabla_{e_1}(j(e_2)) \cdot\overline\varphi + \frac i2 \nabla_{e_1}^\perp (h(e_2))\cdot\overline\varphi}_{\mathcal A^{12}_3}
 \\
&&\underbrace{ - \frac i4 e_2 \cdot j(e_1) \cdot\overline\varphi + \frac i4 j(e_2) \cdot e_1 \cdot\overline \varphi}_{\mathcal A^{12}_4}
\\
&& \underbrace{ - \frac i4 e_2 \cdot h( e_1) \cdot\overline\varphi + \frac i4 h (e_2) \cdot e_1 \cdot\overline\varphi}_{\mathcal A^{12}_5}
\\
&& \underbrace{+ \frac 14 j (e_2) \cdot j (e_1) \cdot \varphi }_{\mathcal A^{12}_6}
+ \underbrace{ \frac 14 h (e_2) \cdot h (e_1 )\cdot \varphi }_{\mathcal A^{12}_7}
\\
&&\underbrace{+ \frac 14 \big(  j (e_2) \cdot h (e_1) \cdot \varphi + h (e_2) \cdot j (e_1) \cdot\varphi\big) }_{\mathcal A^{12}_8}\\
&& \underbrace{- \frac i4 \big(  \eta (e_2)\cdot h (e_1)\cdot \overline\varphi + h (e_2) \cdot \eta(e_1) \cdot\overline \varphi\big)}_{\mathcal A^{12}_9}
\\
&&\underbrace{- \frac i4 \big(  \eta (e_2)\cdot j (e_1)\cdot \overline\varphi + j (e_2) \cdot \eta(e_1) \cdot\overline \varphi\big)}_{\mathcal A^{12}_{10}}.
\end{eqnarray*}

We point out that since $[e_1,e_2]=0$, we have $\nabla_{[e_1, e_2 ]}\varphi =0.$ and some terms are vanishing as shown in the following lemma.
\begin{lemma} We denote by $\mathcal A^{21}_{1}$  (resp. $\mathcal A^{21}_2, \cdots,\mathcal A^{21}_{10})$  the expression $\mathcal A^{12}_1 $ (resp. $\mathcal A^{12}_2, \cdots, \mathcal A^{12}_{10}$) when $e_1$ and $e_2$ are interchanged. 
We have  
\begin{enumerate}
\item 
\begin{eqnarray}
\mathcal A^{12}_2 -  \mathcal A^{21}_{2} = 0
\end{eqnarray}
\item 
\begin{eqnarray}\label{IV}
\mathcal A^{12}_5 -\mathcal  A^{21}_5 =0 
%(e_1, e_2) - IV(e_2, e_1) = 0
\end{eqnarray}
\item 
 \begin{eqnarray}\label{}
\mathcal A^{12}_3 +\mathcal  A^{12}_9 +\mathcal A^{12}_{10} -\mathcal A^{21}_3 - \mathcal A^{21}_9 -\mathcal  A^{21}_{10} = 0
\end{eqnarray}
\item 
\begin{eqnarray}
 \mathcal A^{12}_6- \mathcal  A^{21}_6=-\frac 12 (j_{21})^2e_1\cdot e_2\cdot \varphi
\end{eqnarray}
\item 
\begin{eqnarray}
 \mathcal A^{12}_7 -\mathcal   A^{21}_7=\frac 12 [  h_{21} h_{12}-  h_{11} h_{22} ]\nu_1\cdot \nu_2\cdot\varphi
\end{eqnarray}
\item 
\begin{eqnarray} \mathcal A^{12}_4 - \mathcal  A^{21}_4= i j_{12}  \overline \varphi 
\end{eqnarray}
\item 
\begin{eqnarray}
\mathcal  A^{12}_8- \mathcal A^{21}_8 = & \frac 12 &  \Big ( j_{21}h_{11} e_1 \cdot \nu_1 +  j_{21}h_{12} e_1 \cdot \nu_2 \nonumber \\ && + j_{21}h_{21} e_2 \cdot \nu_1  + j_{12}h_{22} e_2 \cdot \nu_2 \Big) \cdot\varphi
\end{eqnarray}

\item 

\begin{eqnarray}
&& \mathcal A^{12}_1 -\mathcal  A^{21}_1 \nonumber\\ &=& -\frac 12 \sum_{j=1}^2 e_j \cdot \big( (\nabla^{'}_{e_1} B)(e_2, e_j) ) -(\nabla^{'}_{e_2} B)(e_1, e_j) \big)\cdot\varphi \nonumber\\ && + \frac 12 g([S_{\nu_1},  S_{\nu_2}] )(e_1), e_2) \nu_1 \cdot \nu_2 \cdot\varphi \nonumber \\ &&
+\frac{1}{2}\big(|B_{12}|^2- \< B_{11},B_{22}\> \big)e_1\cdot e_2 -\frac 12 e_1\cdot e_2 \cdot \varphi, 
\end{eqnarray}
where $\nabla^{'}$ is the natural connection on $T^*M \otimes T^*M \otimes E$.

%We have \begin{eqnarray*}
%&& \frac 14  \eta (e_2) \cdot e_1\cdot\varphi + \frac 14  e_2 \cdot \eta(e_1) \cdot \varphi \\&& 
%-\frac 14  \eta (e_1) \cdot e_2\cdot\varphi - \frac 14  e_1 \cdot \eta(e_2) \cdot \varphi = 0
%\end{eqnarray*}
\end{enumerate}
\end{lemma}
{\it Proof:} \begin{enumerate}
\item Using the definition of $\eta$, we get  $-\frac 12 B(e_j, X) = e_j\cdot\eta(X) -\eta(X)\cdot e_j$. Hence
\begin{eqnarray*} 
\mathcal A^{12}_2 -  \mathcal A^{21}_2 = -\frac 18 B_{21} + \frac 18 B_{12}= 0.
\end{eqnarray*}
\item \begin{eqnarray}\label{IV}
\mathcal  A^{12}_5 -   \mathcal A^{21}_5  &=& 
  -\frac i4 ( e_2 \cdot h( e_1) - h (e_2) \cdot e_1)\cdot\overline\varphi + \frac i4 (e_1 \cdot h( e_2) - h(e_1) \cdot e_2 ) \cdot \overline \varphi \nonumber \\ 
&=& \frac i4 (2 h_{12}    -2  h_{21} )\cdot\overline \varphi = 0, 
\end{eqnarray}
because $X$ and $h(X)$ are orthogonal for any $X \in \Gamma(TM)$ with respect to the metric $\widetilde{g}$.
 \item First we have 
\begin{eqnarray}\label{II}
&&  \mathcal A^{12}_3  -  \mathcal A^{21}_5 \nonumber \\
&=& \frac i2  \nabla_{e_1} (j (e_2))\cdot \overline \varphi + \frac i2  \nabla^\perp_{e_1} (h (e_2))\cdot \overline \varphi -  \frac i2  \nabla_{e_2} (j (e_1))\cdot \overline \varphi - \frac i2  \nabla^\perp_{e_2} (h (e_1))\cdot \overline \varphi\nonumber \\
& = &
 \frac i2 \Big ( (\nabla_{e_1}j)e_2 \cdot\overline\varphi +   (\nabla_{e_1}h)e_2 \cdot\overline\varphi - 
(\nabla_{e_2}j)e_1 \cdot\overline\varphi - (\nabla_{e_2}h)e_1 \cdot\overline\varphi \Big)
\nonumber\\ & =& 
 \frac i2 \Big (  s(B_{12}) \cdot\overline\varphi + S_{h(e_2)} e_1 \cdot\overline \varphi + t(B_{12})\cdot\overline\varphi -B(e_1, j(e_2)) \cdot\overline\varphi \nonumber \\ && - s(B_{12}) \cdot\overline\varphi -S_{h(e_1)} e_2 \cdot\overline\varphi -t(B_{12}) \cdot\overline\varphi + B(e_2, j (e_1)) \cdot\overline \varphi
\Big ) \nonumber \\ 
&=&  \frac i2 \Big (   S_{h(e_2)} e_1 \cdot\overline \varphi   - S_{h(e_1)} e_2 \cdot\overline\varphi    -B(e_1, j(e_2)) \cdot\overline\varphi  + B(e_2, j (e_1)) \cdot\overline \varphi
\Big ) 
\end{eqnarray}
Moreover, we calculate
\begin{eqnarray}\label{II'}
 -B(e_1, j(e_2)) \cdot\overline\varphi  + B(e_2, j (e_1)) \cdot\overline \varphi \nonumber
&=& - j _{21}  B_{11} \cdot\overline\varphi + j_{12} B_{22} \cdot\overline\varphi
\nonumber\\ &=& 2g(j(e_1), e_2) H \cdot\overline \varphi =  2 j_{12} H \cdot\overline \varphi
\end{eqnarray}
and
\begin{eqnarray}\label{II''}
&&    S_{h(e_2)} e_1 \cdot\overline \varphi   - S_{h(e_1)} e_2 \cdot\overline\varphi   \nonumber \\
&=& -<S_{h(e_1)} e_2, e_1> e_1 \cdot\overline\varphi -<S_{h(e_1)} e_2, e_2> e_2 \cdot\overline\varphi \nonumber\\ 
&&+ <S_{h(e_2)} e_1, e_1> e_1 \cdot\overline\varphi + <S_{h(e_2)} e_1, e_2> e_2 \cdot\overline\varphi 
\nonumber\\ &=& -<B_{21}, h(e_1)> e_1 \cdot\overline \varphi  - <B_{22}, h(e_1)> e_2 \cdot\overline \varphi\nonumber \\
&& + <B_{11}, h(e_2)> e_1 \cdot\overline \varphi +  <B_{12}, h(e_2)> e_2 \cdot\overline \varphi 
%\\&=& <s(II(e_1, e_2)), e_1> e_1 \cdot\overline\varphi -<s(II(e_2, e_2)), e_1> e_2 \cdot\overline\varphi 
%\\&& + <s(II(e_1, e_1)), e_2> e_1 \cdot\overline\varphi  + <s(II(e_1, e_2)), e_2> e_2 %\cdot\overline\varphi 
\end{eqnarray}
In addition we have 
\begin{eqnarray}\label{IIX}
&&  \mathcal A^{12}_9- \mathcal  A^{21}_9 \nonumber\\
&=& \frac i4
\Big (
-e_1 \cdot B_{12} \cdot h(e_1) -e_2 \cdot B_{22} \cdot h(e_1) -h(e_2) \cdot e_1 \cdot B_{11} 
-h(e_2) \cdot e_2 \cdot B_{12}\nonumber \\ 
&& +  e_1 \cdot B_{11} \cdot h(e_2) + e_2 \cdot B_{12} \cdot h(e_2)  
+ h(e_1) \cdot e_1 \cdot B_{12} + h(e_1) \cdot e_2 \cdot B_{22}
\Big ) \cdot \overline \varphi
\nonumber\\ 
&=&   
\frac i4 \Big (  2 <B_{12}, h(e_1)> e_1 +  2 <B_{22}, h(e_1)> e_2 \nonumber \\ && - 2 <B_{12}, h(e_2)> e_2 -  2 <B_{11}, h(e_2)> e_1 \big) \cdot\overline \varphi
 \end{eqnarray}

and 

\begin{eqnarray}\label{IX}
&&  \mathcal A^{12}_{10} - \mathcal  A^{21}_{10} \nonumber \\
&=& - \frac i4 \eta (e_2)\cdot j (e_1)\cdot \overline\varphi - \frac i4  j (e_2) \cdot \eta(e_1) \cdot\overline \varphi + \frac i4 \eta (e_1)\cdot j (e_2)\cdot \overline\varphi + \frac i4  j (e_1) \cdot \eta(e_2) \cdot\overline \varphi \nonumber  \\
&=& 
%\frac i4
%\Big (
%-e_1 \cdot B (e_1, e_2) \cdot j(e_1) \cdot -  e_2 \cdot B (e_2, e_2) \cdot j(e_1) \cdot - j (e_2) \cdot e_1 \cdot B(e_1, e_1) \cdot %\nonumber \\ 
%&& - j (e_2) \cdot e_2 \cdot B(e_1, e_2) \cdot  
% + e_1 \cdot B (e_1, e_1) \cdot j(e_2) \cdot + e_2 \cdot B (e_1, e_2) \cdot j(e_2) \cdot\nonumber\\ 
%&& +j (e_1) \cdot e_1 \cdot B(e_2, e_1)\cdot  + j (e_1) \cdot e_2 \cdot B(e_2, e_2) \cdot 
%\Big ) \overline \varphi
 %\nonumber\\ &=&  
\frac i4
\Big (
e_1 \cdot j (e_1) \cdot  B_{11}  + e_2\cdot j (e_1) \cdot  B_{22} -  j (e_2) \cdot e_1 \cdot  B_{11}- j (e_2) \cdot  e_2 \cdot  B_{12} \nonumber\\ 
&& -  e_1 \cdot j (e_2) \cdot  B_{11} -e_2 \cdot j (e_2) \cdot  B_{12} + j(e_1) \cdot  e_1 \cdot  B_{12} + j(e_1) \cdot e_2 \cdot  B_{22} \Big ) \cdot\overline\varphi    
\nonumber\\
&=& \frac i4 \big (   -2 g(j(e_1), e_2) B_{22}+ 2g(j(e_2), e_1) B_{11}  \big ) \cdot\overline \varphi 
\nonumber \\ &=& - i \  j_{12} H \cdot \overline \varphi.
\end{eqnarray}
Now, replacing  (\ref{II'}) and (\ref{II''}) in  (\ref{II}) and combining together with (\ref{IIX}) and (\ref{IX}), we get the desired result.

\item Since $j$ is antisymmetric, we have $j_{kl} = -j_{lk}$  and so  
\begin{eqnarray}\label{V}
\mathcal  A^{12}_6 - \mathcal  A^{21}_6 
 &=&\frac 14  \big ( j(e_2)\cdot j(e_1) - j(e_1) \cdot j(e_2)\big) \cdot\varphi \nonumber\\ 
%&=& \frac 14 (j_{21} j_{12} e_1\cdot e_2 - j_{12}j_{21} e_2\cdot e_1)\cdot\varphi \nonumber\\
& =& \frac 12 j_{21} j_{12} e_1\cdot e_2  \cdot\varphi = - \frac 12 g(j(e_1), e_2)^2 e_1\cdot e_2 \cdot \varphi
\end{eqnarray}

\item 
\begin{eqnarray*}
 && \mathcal A^{12}_7 -  \mathcal A^{21}_7 \nonumber \\&=&  \frac 14 ( h (e_2) \cdot h (e_1) -  h (e_1) \cdot h (e_2)) \nonumber\\
&=&  \frac 14 \big( - h_{21}h_{11} + h_{21} h_{12} \nu_1 \cdot \nu_2 + h_{22} h_{11} \nu_2 \cdot \nu_1 - h_{22} h_{12} + h_{11} h_{21}\nonumber\\ && - h_{11} h_{22} \nu_1 \cdot \nu_2  - h_{12}h_{21}  -h_{11} h_{22} \nu_1 \cdot \nu_2 - h_{12} h_{21} \nu_2 \cdot \nu_1 + h_{12} h_{22}\big) \cdot\varphi 
\nonumber\\ &=& \frac 12 (h_{21} h_{12}  -  h_{11} h_{22}) \nu_1 \cdot \nu_2 \cdot \varphi
\end{eqnarray*}

\item 

\begin{eqnarray}\label{III}
\mathcal  A^{12}_4 -  \mathcal A^{21}_4 &=&- \frac i4 \big ( e_2 \cdot j(e_1) - j(e_2) \cdot e_1  -e_1 \cdot j(e_2) +j(e_1) \cdot e_2\big) \cdot\overline\varphi \nonumber \\
&=& - \frac i4 \big ( - j_{12}+ j_{21} +  j_{21} - j_{12}\big) \cdot\overline\varphi  = i \ j_{12} \overline \varphi.
%&=& \frac i4 (2 g(e_2, j(e_1))   -2 g(e_1, j(e_2))  ) \cdot\overline \varphi\nonumber\\
%&=& i g(e_2, j (e_1)) \overline \varphi  
\end{eqnarray}

\item 
We have

\begin{eqnarray}\label{VII}
&& \mathcal  A^{12}_8 -  \mathcal A^{21}_8\nonumber \\
&=& \frac 14 \big(  j (e_2) \cdot h (e_1)  + h (e_2) \cdot j (e_1)  -   j (e_1) \cdot h (e_2)  - h (e_1) \cdot j (e_2)  \big) \cdot \varphi \nonumber\\
%&=&  
%\frac 14 \big(  j_{21} e_1 \cdot (h_{11} \nu_1 + h_{12} \nu_2) + j_{12}(h_{21} \nu_1 + h_{22} \nu_2) e_2 \nonumber\\ && - j_{12} e_2 %\cdot(h_{21} \nu_1 + h_{22} \nu_2) - j_{21} (h_{11} \nu_1 + h_{12} \nu_2) e_1 \big) \cdot \varphi
%\nonumber\\ 
&=& \frac 12 \big ( j_{21}h_{11} e_1 \cdot \nu_1 + j_{21}h_{12} e_1 \cdot \nu_2 + j_{21}h_{21} e_2 \cdot \nu_1
+ j_{12}h_{22} e_2 \cdot \nu_2 \big) \cdot\varphi
\end{eqnarray}

Finally,  we recall here that Lemma 3.3 of \cite{BLR} says that first
$$d\eta(X,Y)=-\frac 12 \sum_{j=1}^2 e_j \cdot \big( (\nabla^{'}_{X} B)(Y, e_j) ) -(\nabla^{'}_{Y} B)(X, e_j) \big),$$
and second 
\begin{eqnarray*}
\eta(e_2)\cdot\eta(e_1)-\eta(e_1)\cdot\eta(e_2)&=&\frac{1}{2}\big(|B_{12}|^2-\left<B_{11},B_{22}\right>\big)e_1\cdot e_2\\
&&+\frac{1}{2}\left<\left(S_{\nu_1}\circ S_{\nu_2}-S_{\nu_2}\circ S_{\nu_1}\right)(e_1),e_2\right>\nu_1\cdot \nu_2.
\end{eqnarray*}
Moreover,  since $d\eta(e_1,e_2)=\nabla_{e_1}(\eta(e_2))-\nabla_{e_2}(\eta(e_1))$ , we deduce immediately from the definition of $\mathcal A^{12}_1$ and  $\mathcal  A^{21}_1$ and the two above identities the desired relation.
\end{enumerate}
\hfill$\square$\\ 
Now, we have all the ingredients to prove Theorems \ref{thm1} and \ref{thm2}. 
\section{Lagrangian case, proof of Theorem \ref{thm2}}
First and from Sections \ref{sec3} and \ref{sec4}, assertion (2) of Theorem \ref{thm2} implies assertion (1). Assume now that assertion (1) is satisfied.  Since, $j=t=0$, we have 
\begin{eqnarray}
\mathcal{R}_{e_1, e_2} \varphi &=& \frac 12 K_M e_1 \cdot e_2 \cdot \varphi - \frac 12 K_E \nu_1 \cdot \nu_2 \cdot \varphi + \frac 12 F^{M +E} (e_1, e_2) \varphi,
\end{eqnarray}
with $F^{M +E} (e_1, e_2)  =0$  because $j=0$.
On the other hand, we have 

\begin{eqnarray}
 \mathcal{R}_{e_1, e_2} \varphi 
&=& - \frac 12 \sum_{j=1}^2 e_j \cdot( (\nabla^{'}_{e_1} B) (e_2, e_j)     -  (\nabla^{'}_{e_2} B) (e_1, e_j)    ) \cdot\varphi \nonumber \\
&& + \frac 12 (\vert B_{12}\vert^2 - <B_{11}, B_{22}>) e_1 \cdot e_2 \cdot \varphi \nonumber\\
&& + \frac 12 <[S_{\nu_1}, S_{\nu_2}](e_1), e_2> \nu_1 \cdot \nu_2 \cdot\varphi \nonumber \\
&&-\frac 12 e_1 \cdot e_2 \cdot \varphi + \frac 12  (h_{21}h_{12} - h_{11}h_{22}) \nu_1 \cdot\nu_2 \cdot\varphi
\end{eqnarray}

We get finally that $T \cdot \varphi = 0$, where $T \in (\Lambda^2 M \otimes 1 \oplus TM \otimes E \oplus 1 \otimes \Lambda^2 E)$ is given by
\begin{eqnarray*}
T&=&\frac12(<B_{11}, B_{22}> - \vert B_{12}\vert^2 +1-K_M)e_1\wedge e_2\\
&&+\frac12(h_{22}h_{11}- h_{21} h_{12} -<[S_{\nu_1}, S_{\nu_2}](e_1), e_2>  -K_E)\nu_1\wedge\nu_2\\
&&-\frac12\sum_{j=1}^2 e_j \wedge( (\nabla^{'}_{e_1} B) (e_2, e_j)     -  (\nabla^{'}_{e_2} B) (e_1, e_j)    ) .
\end{eqnarray*} 
Now, we recall that Lemma 3.4 of \cite{BLR} ensures that if $T$ is a two form and $\varphi$ a spinor so that  $\varphi^+$ and $\varphi^-$ never vanish and $T\cdot\varphi=0$, then $T=0$. Note that the hypothesis that both $\varphi^+$ and $\varphi^-$ do not vanish is crucial. Here, the conclusion $T=0$ reduces to the following identities

$$K_M = <B_{11}, B_{22}> - \vert B_{12}\vert^2 +1,$$
$$ K_E = - <[S_{\nu_1}, S_{\nu_2}](e_1), e_2> - (h_{21} h_{12} - h_{22}h_{11} ),$$
$$(\nabla^{'}_{e_1} B) (e_2, e_j)     -  (\nabla^{'}_{e_2} B) (e_1, e_j)   = 0, $$
which are Gauss, Ricci and Codazzi equations for a Lagrangian surface in $\mathbb C P^2$ and so the conditions \eqref{compatibilitylagrangian} are fulfilled. Hence, by Proposition \ref{propimmersioncp2}, we conclude that there exists a Lagrangian isomertic immersion from $(M,g)$ into $\CC P^2$ with $E$ as normal bundle and $B$ as second fundamental form. This proves that assertion (1) of Theorem \ref{thm2} implies assertion (2). Theorem \ref{thm2} is proved.

\section{Complex case, proof of Theorem \ref{thm1}}
Again, assertion (2) of Theorem \ref{thm1} implies assertion (1) by the discussions of Sections \ref{sec3} and \ref{sec4}. Assume now that assertion (1) is satisfied.  We have  $s =0$, $ h = 0$  so $F^{M+E} (e_1, e_2) = -2i$. 
We take $j (e_1) = e_2$ and $t \nu_1 = \nu_2$, \emph{i.e.} $g(j (e_1), e_2) = g (t \nu_1, \nu_2 )=1.$  We calculate and we get 
\begin{eqnarray}
\mathcal{R}_{e_1, e_2} \varphi &=&- \frac 12 K_M e_1 \cdot e_2 \cdot \varphi - \frac 12 K_N \nu_1 \cdot \nu_2 \cdot \varphi + \frac 12 F^{M +E} (e_1, e_2) \varphi \nonumber \\ &=&
- \frac 12 K_M e_1 \cdot e_2 \cdot \varphi - \frac 12 K_N \nu_1 \cdot \nu_2 \cdot \varphi -i  \varphi  \nonumber \\ &=& 
\overline T \cdot \varphi - i \varphi,
\end{eqnarray}
where $\overline{T}$ is the $2$-form defined by
\begin{eqnarray*}
\overline{T}&=&- \frac 12 K_M e_1 \wedge e_2 - \frac 12 K_N \nu_1 \wedge\nu_2
\end{eqnarray*}
On the other hand, we have
\begin{eqnarray}
\mathcal{R}_{e_1, e_2} \varphi &=& - e_1 \cdot e_2 \cdot\varphi + i \overline \varphi 
\nonumber \\
&& - \frac 12 \sum_{j=1}^2 e_j \cdot( (\nabla'_{e_1} B) (e_2, e_j)     -  (\nabla'_{e_2} B) (e_1, e_j)    ) \cdot\varphi \nonumber \\
&& + \frac 12 (\vert B_{12}\vert^2 - <B_{11}, B_{22}>) e_1 \cdot e_2 \cdot \varphi \nonumber\\
&& + \frac 12 <[S_{\nu_1}, S_{\nu_2}](e_1), e_2> \nu_1 \cdot \nu_2 \cdot\varphi \nonumber \\ 
&& = \widetilde T \cdot \varphi + i \overline \varphi,
\end{eqnarray}
where $\widetilde{T}$ is the $2$-form defined by
\begin{eqnarray*}
\widetilde{T}&=&\frac12(<B_{11}, B_{22}> - \vert B_{12}\vert^2 +1-K_M)e_1\wedge e_2\\
&&+\frac12(h_{22}h_{11}- h_{21} h_{12} -<[S_{\nu_1}, S_{\nu_2}](e_1), e_2>  -K_E)\nu_1\wedge\nu_2\\
&&-\frac12\sum_{j=1}^2 e_j \wedge( (\nabla^{'}_{e_1} B) (e_2, e_j)     -  (\nabla^{'}_{e_2} B) (e_1, e_j)    ) .
\end{eqnarray*} 
Together, it gives $\overline T \cdot \varphi - \widetilde T \cdot\overline \varphi - i \varphi - i \overline \varphi = 0$, which means that $\mathcal T \cdot \varphi - i \varphi - i \overline \varphi = 0$, where the $2$-form $\mathcal T=\overline{T}-\widetilde{T}$ is  given by 
\begin{eqnarray}\label{TT}
\mathcal T &=& - \frac 12 K_M e_1 \wedge e_2 \cdot \varphi - \frac 12 K_N \nu_1 \wedge \nu_2 \cdot \varphi\nonumber + e_1\wedge e_2 \cdot\varphi \\ 
&&  -\frac 12 (\vert B_{12}\vert^2 - <B_{11}, B_{22}>) e_1 \wedge e_2 \cdot \varphi \nonumber\\
&& - \frac 12 <[S_{\nu_1}, S_{\nu_2}](e_1), e_2> \nu_1 \wedge \nu_2 \cdot\varphi \nonumber \\ 
&& +\frac 12 \sum_{j=1}^2 e_j \wedge( (\nabla'_{e_1} B) (e_2, e_j)     -  (\nabla'_{e_2} B) (e_1, e_j)    ) \cdot\varphi
\end{eqnarray}
We give now the following Lemma:
\begin{lemma}
Let $\mathcal T$ be a $2$ form, \emph{i.e.} $\mathcal T \in \Lambda^2 M \otimes 1 \oplus \Lambda^1 M \otimes \Lambda^1 E \oplus 1 \otimes \Lambda^2 E$ and $\varphi\in\Sigma$ so that both $\varphi^+$ and $\varphi^-$ never vanish. Assume that 
$$\mathcal T \cdot \varphi - i \varphi - i \overline \varphi = 0, $$
and write $\mathcal T = T^t e_1 \wedge e_2  + T^n \nu_1 \wedge \nu_2 + T^m,$ where $T^m \in \Lambda^1 M \otimes \Lambda^1E$. 
Then, 
$$T^t = -1, T^n = 0 \ \ \text{and}\ \  T^m =  0 .$$
\end{lemma}
{\bf Proof.} 
Let $\varphi = \varphi^+ + \varphi^-,$ with
$$\varphi^+ =\varphi^{++} + \varphi^{--},$$
$$\varphi^- = \varphi^{-+} + \varphi^{+-},$$
a solution of \eqref{partspinor} with $h=0$. This means that 
 $$\nabla_X \varphi^{++} = -\frac 12 X\cdot \varphi^{-+} - \frac i2 j(X)\cdot \varphi^{-+}$$
  $$\nabla_X \varphi^{+-} = -\frac 12 X\cdot \varphi^{--} + \frac i2 j(X)\cdot \varphi^{--}$$
 $$\nabla_X \varphi^{-+} = -\frac 12 X\cdot \varphi^{++} +\frac i2 j(X)\cdot \varphi^{++}$$
  $$\nabla_X \varphi^{--} = -\frac 12 X\cdot \varphi^{+-} - \frac i2 j(X)\cdot \varphi^{+-}.$$
  For a sake of simplicity, and without lost of generality, we can restrict only $\varphi^+=\varphi^{++}$ and $\varphi^-=\varphi^{-+}$ which which have no zeros by assumption. The equation $$\mathcal T \cdot \varphi - i \varphi - i \overline \varphi = 0,$$ 
  becomes $$T^t e_1\cdot e_2\cdot (\varphi^{++} + \varphi^{-+}) + (T^n+1) \nu_1\cdot \nu_2 \cdot (\varphi^{++} + \varphi^{-+}) + T^m \cdot (\varphi^{++} + \varphi^{-+}) = i \overline \varphi = i (\varphi^{++} - \varphi^{-+})$$
 Taking  the scalar product with $\varphi^{++}$ then with $\varphi^{-+}$, we get
$$T^t +T^n +1 = -1,$$
$$-T^t + T^n+1 = 1,$$
  which gives $T^n =-1$, $T^t = -1$ and $T^m = 0$. One can uses (\ref{TT}) to get Gauss, Codazzi and  Ricci equations and so the conditions \eqref{compatibilitycomplex} are fulfilled. There are exactly the conditions of a complex immersion. Hence, by Proposition \ref{propimmersioncp2}, we conclude that there exists a complex isometric immersion from $(M,g)$ into $\CC P^2$ with $E$ as normal bundle and $B$ as second fundamental form. As for the Lagrangian case, this proves that assertion (2) of Theorem \ref{thm1} implies assertion (1). 

%\textcolor{red}{Il est sans aucun doute possible de d\'emontrer une version avec l'op\'erateur de Dirac et une condition sur la d\'eriv\'ee de la norme comme je l'ai fait pour les produits dans \cite{Roth4}. \'A d\'efaut de la faire en d\'etails, on peut le citer comme corollaire et donner quelques arguments de preuve et se r\'ef\'erer \`a \cite{Roth4} dans une derni\`ere section. Qu'en penses-tu?\\ \\
%Je viens de regarder rapidement, je pense que \c{c}a marche sans probl\`eme avec une toute petite adaptation. Le probl\`eme est que la preuve compl\`ete est assez longue, 6 pages environ dans mon papier sur les produits. Soit on la met en entier, mais le papier devient long (plus de 25 pages...) soit on \'enonce en disant que c'est similaire sans donner la preuve. Je suis plut\^ot pour la seconde option car sinon le papier risque d'\^etre trop long. Au pire, le moment venu, le referee peut nous demander d'ajouter la preuve.}

\section{The Dirac equation}
Let $\varphi$ be a spinor field satisfying Equation \eqref{partspinor}, then it satisfies the following Dirac equation
\begin{equation}\label{eqdirac}
D\varphi=\vec{H}\cdot\varphi-\varphi +\frac i2\beta\cdot\overline{\varphi},
\end{equation}
where $\beta$ is defined by $ \displaystyle\beta=\sum_{i=1,2}e_i\cdot h(e_i)=\sum_{i,j=1}^2h_{ij}e_i\cdot\nu_j$, where $\{e_1,e_2\}$ and $\{\nu_1,\nu_2\}$ are respectively orthonormal frames of $TM$ and $E$ and $h_{ij}=\langle h(e_i),\nu_j\rangle$.\\
As in \cite{BLR} and \cite{Roth4}, we will show that this equation with an appropiate condition on the norm of both $\varphi^+$ and $\varphi^-$ is equivalent to Equation \eqref{partspinor}, where the tensor $B$ is expressed in terms of the spinor field $\varphi$ and such that $\tr(B)=2\vec{H}$. Moreover, from 
$$
\nabla_X\varphi^{\pm}=-\frac12\eta(X)\cdot\varphi^{\pm}-\frac12 X\cdot\varphi^{\mp}\mp\frac i2j(X)\varphi^{\mp}\mp\frac i2h(X)\varphi^{\mp},
$$
we deduce that
\begin{equation}\label{normepm}
X(|\varphi^{\pm}|^2)=\pre \left\langle -\frac12 X\cdot\varphi^{\mp}\mp\frac i2j(X)\cdot \varphi^{\mp}\mp\frac i2h(X)\cdot\varphi^{\mp},\varphi^{\pm}\right\rangle.
\end{equation}
Now, let $\varphi$ a spinor field solution of the Dirac equation \eqref{eqdirac} with $\varphi^+$ and $\varphi^-$ nowhere vanishing and satisfying the norm condition \eqref{normepm}, we set for any vector fields $X$ and $Y$ tangent to $M$ and $\xi\in \Gamma(E)$
\begin{eqnarray}\label{defB}
&& \left<B(X,Y),\xi\right> \nonumber \\&=&\frac{1}{|\varphi^+|^2}\pre\left<X\cdot\nabla_Y\varphi^+-\frac12 \left(X+ij(X)+ih(X)\right)\cdot Y\cdot\varphi^{-},\xi\cdot\varphi^+\right> \nonumber\\
&&+\frac{1}{|\varphi^-|^2}\pre\left<X\cdot\nabla_Y\varphi^--\frac12 \left(X-ij(X)-ih(X)\right)\cdot Y\cdot\varphi^-,\xi\cdot\varphi^+\right>
\end{eqnarray}

Then, we have the following
\begin{prop}\label{prop1}
Let $\varphi\in\Gamma(\Sigma)$ satisfying the Dirac equation 
$$D\varphi=\vec{H}\cdot\varphi-\varphi+\frac i2 \beta\cdot\overline{\varphi}$$
such that 
$$
X(|\varphi^{\pm}|^2)=\pre \left\langle -\frac12 X\cdot\varphi^{\mp}\mp\frac i2j(X)\cdot\varphi^{\mp}\mp\frac i2h(X)\cdot\varphi^{\mp},\varphi^{\pm}\right\rangle,
$$
then $\varphi$ is solution of Equation \eqref{partspinor}
$$\nabla_X\varphi=-\frac 12\eta(X)\cdot\varphi-\frac{1}{2}X\cdot\varphi+\frac{i}{2}j(X)\cdot\overline{\varphi}+\frac{i}{2}h(X)\cdot\overline{\varphi},$$
where $\eta$ is defined by $\displaystyle\eta(X)=\sum_{j=1}^2e_j\cdot B(e_j,X)$. Moreover, $B$ is symmetric.
\end{prop}
The proof of this proposition will not be given, since it is completely similar to the case of Riemannian products \cite[Proposition 4.1]{Roth4}. Now, combining this proposition with Theorems \ref{thm1} and \ref{thm2}, we get the following corollaries. We have this first one for complex immersions of surfaces.
\begin{cor}\label{cor1}
Let $(M^2,g)$ be an oriented Riemannian surface and $E$ an oriented vector bundle of rank $2$ over $M$ with scalar product $<\cdot,\cdot>_E$ and compatible connection $\nabla^E$. We denote by $\Sigma=\Sigma M\otimes\Sigma E$ the twisted spinor bundle. Let $j$ be a complex structure on $M$ and $t$ a complex structure on $E$. Let $\vec{H}$ be a section of $E$. Then, the following two statements are equivalent
\begin{enumerate}

\item There exists a ${\rm Spin}^c$  structure on $\Sigma M\otimes\Sigma E$ whose auxiliary line bundle's curvature is given by $F^{M +E} (e_1, e_2)=0$ and a spinor field $\varphi$ in $\Sigma$ solution of the Dirac equation
$$D\varphi=\vec{H}\cdot\varphi-\varphi$$
such that $\varphi^+$ and $\varphi^-$ never vanish satisfying the norm condition
$$X(|\varphi^{\pm}|^2)=\pre \left\langle -\frac12 X\cdot\varphi^{\mp}\mp\frac i2j(X)\cdot \varphi^{\mp}\varphi^{\pm}\right\rangle$$
 and such that the maps $j$, $t$ and the tensor $B$ defined by \eqref{defB} satisfies $t(B(X,Y))=B(X,j(Y))$ for all $X,Y\in\Gamma(TM)$.
 \item There exists an isometric {\bf complex} immersion of $(M^2,g)$ into $\CC P^2$ with $E$ as normal bundle and mean curvature $\vec{H}$ such that over $M$ the complex structure of $\CC P^2$ is given by $j$ and $t$ (in the sense of Proposition \ref{propimmersioncp2}).
\end{enumerate}
\end{cor}

We have this second corollary for Lagrangian surfaces.

\begin{cor}\label{cor2}
Let $(M^2,g)$ be an oriented Riemannian surface and $E$ an oriented vector bundle of rank $2$ over $M$ with scalar product $<\cdot,\cdot>_E$ and compatible connection $\nabla^E$. We denote by $\Sigma=\Sigma M\otimes\Sigma E$ the twisted spinor bundle. Let $B:TM\times TM\lgra E$ a bilinear symmetric map, $h:TM\lgra E$ and $s: E\lgra TM$ the dual map of $h$. Assume that the maps $h$, $s$ are parallel and satisfy $h\circ s=-\id_{E}$. Let $\vec{H}$ be a section of $E$. Then, the following two statements are equivalent
\begin{enumerate}

\item There exists a ${\rm Spin}^c$  structure on $\Sigma M\otimes\Sigma E$ whose auxiliary line bundle's curvature is given by $F^{M +E} (e_1, e_2)=-2i$ and a spinor field $\varphi$ in $\Sigma$ solution of the Dirac equation
$$D\varphi=\vec{H}\cdot\varphi-\varphi+\frac i2\beta\cdot\overline{\varphi}$$
\big($\beta$ is the 2-form defined by $ \displaystyle\beta=\sum_{i=1,2}e_i\cdot h(e_i)$\big)  
such that $\varphi^+$ and $\varphi^-$ never vanish  satisfying the norm condition
$$X(|\varphi^{\pm}|^2)=\pre \left\langle -\frac12 X\cdot\varphi^{\mp}\mp\frac i2h(X)\cdot \varphi^{\mp}\varphi^{\pm}\right\rangle$$
 and such that  the tensor $B$ defined by \eqref{defB} satisfies $A_{h(Y)}X+s(B(X,Y))=0$, for all $X\in TM$, where $A_{\nu}:TM\lgra TM$ if defined by $g(A_{\nu}X,Y)=\langle B(X,Y),\nu\rangle_E$ for all $X,Y\in \Gamma(TM)$ and $\nu\in \Gamma(E)$.
 \item There exists an isometric {\bf Lagrangian} immersion of $(M^2,g)$ into $\CC P^2$ with $E$ as normal bundle and mean curvature $\vec{H}$ such that over $M$ the complex structure of $\CC P^2$ is given by $h$ and $s$ (in the sense of Proposition \ref{propimmersioncp2}). 
\end{enumerate}
\end{cor}
{\bf Acknowledgment.}  
The first named author would like to thank the University of Paris-Est, Marne La Vall\'ee  for its support and hospitality. The authors warmly thank the referee for his/her remraks that allow to highly improve the present paper. Both authors are also grateful to Mihaela Pilca for helpful discussions about K\"ahlerian Killing spinors.

\end{document}